\documentclass[graybox]{svmult}

\usepackage{mathptmx}       
\usepackage{helvet}         
\usepackage{courier}        
\usepackage{type1cm}        
\usepackage{makeidx}         
\usepackage{graphicx}        
\usepackage{multicol}        
\usepackage[bottom]{footmisc}

\usepackage{tabularx, booktabs}
\usepackage{tikz}
\usepackage{float}
\usepackage{array}
\usepackage{caption}
\usepackage{amssymb,amsmath}
\usetikzlibrary{arrows}

\newcommand\eps{\varepsilon}
\newcommand\C{\mathbb C}

\renewcommand\P{\mathbb P}

\newcommand\Z{\mathbb Z}
\newcommand\calo{{\mathcal O}}

\newcommand\calh{{\mathcal H}}

\newcommand\cala{{\mathcal A}}

\def\endproof{\hspace*{\fill}\endproofsymbol\endtrivlist}

\def\endproofsymbol{\frame{\rule[0pt]{0pt}{6pt}\rule[0pt]{6pt}{0pt}}}

\newcommand\numberthis{\addtocounter{equation}{1}\tag{\theequation}}

\newtheorem{observation}{Observation}

\newcommand\newop[2]{\def#1{\mathop{\rm #2}\nolimits}}
\newop\Ass{Ass}
\newop\Zeroes{Zeroes}
\newop\GL{GL}
\newop\HF{HF}
\newop\Der{Der}
\newop\reg{reg}

\makeindex             

\begin{document}

\title*{Fermat-type arrangements}
\author{Justyna Szpond}
\institute{Justyna Szpond \at Pedagogical University of Cracow,
Podchor\c a\.zych 2, PL-30-084 Cracow, Poland
\email{szpond@gmail.com}}
%
%
\maketitle

\abstract{The purpose of this work is to collect in one place available information on line arrangements
   known in the literature as braid, monomial, Ceva or Fermat arrangement. They have been studied for
   a long time and appeared recently in connection with highly interesting problems, namely:
   the containment problem between symbolic and ordinary powers of ideals and the existence
   of unexpected hypersurfaces. We also study also derived configurations of points
   (or more general: linear flats) which arise by intersecting hyperplanes in Fermat
   arrangements or by taking duals of these hyperplanes. Furthermore we discuss briefly higher dimensional generalizations
   and present results arising by applying this approach to problems mentioned above. Some of our results
   are original and appear for the first time in print.}

\section{Introduction}
\label{JSsec:1}
   Fermat arrangements of lines, as reflection arrangements appear (under the name of Ceva arrangements)
   in Hirzebruch's work \cite{Hir83}. Hirzebruch's interest in them was motivated by seeking ways
   to construct surfaces of general type which are ball quotients. This treatment has been
   considerably extended in the book by Barthel, Hirzebruch and H\"ofer \cite{BHH87}.
   For a recent update on relations between ball quotients and line arrangements we refer to
   Tretkoff's book \cite{Tretkoff16}.

   It seems that the name Fermat arrangements first appears in PhD thesis of Urzua \cite{UrzuaPhD}.
   They came into focus in connection with examples of the non-containment between the third
   symbolic and the second ordinary power of an ideal of a set of points in $\P^2$, which appeared
   first in the work of Dumnicki, Szemberg and Tutaj-Gasi\'nska \cite{DST13} and were
   considerable generalized by Harbourne and Seceleanu \cite{HarSec15}.
   Subsequently asymptotic invariants of associated ideals of points have been computed
   in \cite{DHNSST15}.
   Nagel and Seceleanu
   provided in \cite{NagSec16} detailed description of Rees algebras of ideals of points
   derived from Fermat arrangements of lines.

   The idea of studying generalizations of Fermat line arrangements to Fermat-type arrangements of hyperplanes
   in higher dimensional projective spaces appeared during the workshop
   ''Ordinary and Symbolic Powers of Ideals'' held in Oaxaka in May 2017, particularly
   during discussions with Juan Migliore and Uwe Nagel which I enjoyed so much.
   My joint papers with Malara \cite{MalSzp18}, \cite{MalSzp18Moieciu} show how
   configurations of codimension two flats derived from Fermat-type arrangements provide
   additional non-containment results.

   Additional interest in examples studied here comes from yet another direction.
   The Bounded Negativity Conjecture predicts that on any complex algebraic surface
   the self-intersection numbers of reduced and irreducible curves are bounded from below.
   This Conjecture is one of the central and most difficult problems in the theory
   of algebraic surfaces. Seminal work of Bauer, Harbourne, Knutsen, K\"uronya,
   M\"uller-Stach, Roulleau, and Szemberg \cite{BNC} renewed interest in this conjecture.
   The subsequent article by Bauer, Di Rocco, Harbourne, Huizenga, Lundman,
   Pokora,
   and Szemberg \cite{BNCarrangements15} revealed a link between hyperplane
   arrangements and the Bounded Negativity Conjecture. Fermat-type arrangements
   served there as examples with extremal Harbourne constants.

   My initial idea was to treat here all this interesting developments.
   I realized soon that this would exceed by far the scope of a conference proceedings
   article. Therefor I decided to focus on the latest theory where Fermat-type
   arrangements and derived configurations seem to play a prominent role.
   I mean here the theory of unexpected hypersurfaces initiated in the edge-cutting
   work of Cook II, Harbourne, Migliore, and Nagel \cite{CHMN} developing
   very rapidly. Thus after collecting some general facts about reflection
   arrangements I pass directly to Fermat-type arrangements and unexpected
   hypersurfaces. Some of results presented there are new and hopefully
   so surprising, that they will ignite new pathes of research.

   I failed in writing a comprehensive survey on Fermat-type arrangements
   and their appearances in commutative algebra and algebraic geometry. Maybe some other time.
\subsection{Notation}
   We adopt the combinatorial convention and write
   $$[t]=\left\{1,\ldots,t\right\}.$$
   Even though many arguments are valid over an arbitrary field containing enough roots of unity,
   in order to avoid additional assumptions at some places and to streamline the discussion,
   we make a general assumption of working solely over the field of complex numbers.

\section{Reflection groups and arrangements}
\label{JSsec:2}
   In this section we recall briefly some general facts about arrangements associated to finite reflection groups.
   This will serve as a motivation for one of possible generalizations of Fermat arrangements of lines.

   Let $V$ be the unitary space $\C^{N+1}$ with standard basis $e_0,\ldots,e_N$ and unitary product $<\cdot,\cdot>$.
   Let $R=\C[x_0,\ldots,x_N]$ be the ring of polynomials. Let $G$ be a finite
   group of linear automorphisms of $V$. Then $G$ acts on $R$ by
   $$(g\cdot f)(x)=f(g^{-1}(x))$$
   for all $x\in V$, $f\in R$ and $g\in G$.
\begin{definition}[Imprimitive group of automorphisms]
   A group $G$ of unitary automorphisms of $V$ is called \emph{imprimitive},
   if $V$ is a direct sum
   \begin{equation}\label{eq: imprimitive}
      V=V_1\oplus\ldots\oplus V_t
   \end{equation}
   of non-trivial proper linear subspaces $V_i$ of $V$ such that the set
   $\left\{V_i, i\in[t]\right\}$ is invariant under $G$. This set is called
   a \emph{system of imprimitivity} for $G$.
\end{definition}
   We are mainly interested here in reflection groups.
\begin{definition}[Reflections and reflection groups]
   A linear automorphism $s$ of $V$ of finite order is a \emph{reflection} in $V$
   if it has exactly $N$ eigenvalues equal to $1$.\\
   A \emph{reflection group} in $V$ is a group generated by reflections in $V$.
\end{definition}
   Equivalently, a reflection is a non-trivial linear automorphism of finite order that fixes a hyperplane.
   We call this hyperplane a \emph{reflection hyperplane}.
\begin{remark}
   Any reflection $s$ in $V$ of order $d\geq 2$ has the form
   $s=s_{a,\eps}$ with
   $$s_{a,\eps}(x)=x-(1-\eps)\frac{<x,a>}{<a,a>}a$$
   for some vector $a\in V$ and a primitive root of unity $\eps$ of order $d$.
\end{remark}
\begin{remark}
   If $V$ does not admit any decomposition as in \eqref{eq: imprimitive}, then
   $G$ is called primitive.
\end{remark}
\begin{definition}[Irreducible group]
   A group $G\subset\GL(V)$ is called \emph{irreducible} if there is no non-trivial
   invariant proper subspace of $V$ invariant under $G$.
\end{definition}
   The following example is important for our considerations.
\begin{example}[Monomial gropus]
   Let $\Pi_{N+1}\subset\GL(V)$ be the group of all $(N+1)\times (N+1)$ permutation matrices.
   It is of course isomorphic with the permutation group $S_{N+1}$ of $(N+1)$ elements.
   Let $n\geq 2$ and $p\geq 1$ be integers with $p|n$ and let $A(n,p,N+1)$ be the group
   of $(N+1)\times (N+1)$ diagonal matrices $A=\left(a_{ij}\right)_{i,j\in[N+1]}$
   with $a_{ij}=\eps^{\alpha_i}\delta_{ij}$, where
   $\eps$ is a primitive root of unity of order $n$,
   $\alpha_i\in[n]$ and such the product
   $$\det(A)=\prod_{i\in[N+1]}a_{ii}$$
   is a power of $\eps^p$.
   Let $G(n,p,N+1)$ be the semi-direct product of $A(n,p,N+1)$ and $\Pi_{N+1}$.
   Then $G(n,p,N+1)$ is an irreducible unitary reflection group.
   The generators are matrices in $A(n,p,N+1)$ whose all but exactly one diagonal entries
   are $1$ and products of these matrices with permutation matrices.
\end{example}
   Expressed somewhat simpler, $G(n,p,N+1)$ consists of all monomial $(N+1)\times (N+1)$
   matrices whose non-zero elements are roots of unity of order $n$ and the product of all of these entries
   is a root of unity of order $n/p$. Recall that a matrix is called \emph{monomial} if in any row and any column
   there is exactly one non-zero element. In particular permutation matrices are monomial.\\
   Allowing $n=1$ we could identify $\Pi_{N+1}$ with the group $G(1,1,N+1)$. Note that
   this group is reducible.

   There are inclusions
   $$G(1,1,N+1)\subset G(n,p,N+1)\subset G(n,1,N+1).$$
   The group $G(n,1,N+1)$ is called the \emph{full monomial group}, see \cite[Example 6.29]{OrlTer}
   or \cite[Example 2.23]{dim}.
\begin{remark}\label{rmk: GnpN hyperplanes}
   The reflection hyperplanes for $G(n,p,N+1)$ with $p<n$ are of the form
   $$x_i=\eps^\alpha x_j$$
   for $0\leq i<j\leq N$ and $\alpha\in[n]$
   or the coordinate hyperplanes
   $$x_i=0$$
   for $i\in\left\{0,\ldots,N\right\}$.\\
   For $G(n,n,N+1)$ only the hyperplanes of the first kind are reflection hyperplanes.
\end{remark}

   Finite complex reflection groups were classified by Shephard and Todd \cite{SheTod54}.
   In particular they showed the following result \cite[Section 2]{SheTod54}, see also \cite[Theorem 2.4]{Coh76}
   for an alternative proof.
\begin{theorem}[Imprimitive reflection groups]
   The only (up to conjugation) irreducible imprimitive unitary groups generated by reflections are
   $$G(n,p,N+1)$$
   for $n\geq 2$, $N\geq 1$ and $p|n$ with exception of $G(2,2,2)$ which is reducible.
\end{theorem}
\begin{remark}
   The groups denoted here by $G(n,p,N+1)$ are usually in the literature denoted by $G(m,p,n)$.
   We thought it less confusing to use  right away the notation which applies in other parts of this article.
\end{remark}
\begin{definition}[Semi-invariant polynomials]
   A polynomial $f\in R$ is called \emph{semi-invariant} with respect to $G$,
   if there exists a linear character $\eta$ of $G$ such that
   $$g\cdot f=\eta(g)\cdot f$$
   for all $g\in G$.
\end{definition}
\begin{remark}
   It is well known (see e.g. \cite[Proposition 2.2]{Coh76}
   that if $G$ is an imprimitive and irreducible finite reflection group and $N$ is at least $1$ (i.e. $\dim(V)\geq 2$),
   then $\dim(V_i)=1$ for all $i\in[t]$, hence $t=N+1$ and there are mutually
   distinct homogeneous linear polynomials $\ell_1,\ldots,\ell_t$ such that
   their product $\ell_1\cdot\ldots\cdot\ell_t$ is a homogeneous, semi-invariant with respect to $G$,
   polynomial of degree $N+1$ in $R$.
\end{remark}
\begin{definition}[Reflection arrangement]
   To any finite reflection group $G$ one associates a hyperplane arrangement $\calh(G)$,
   which consists of reflection hyperplanes defined by elements in $G$.
\end{definition}
\begin{example}[Braid arrangement]\label{ex: braid}
   The group $G(1,1,N+1)$ is just a representation of the symmetry group $S_{N+1}$ acting on $V=\C^{N+1}$
   by permuting coordinates. The reflecting hyperplanes are given by equations
   $$x_i=x_j$$
   for $0\leq i<j\leq N$. This arrangement can be studied also projectively. The hyperplanes
   in $\P^N$ are defined by linear factors of the semi-invariant polynomial
   $$F_{N,1}=\prod_{0\leq i<j\leq N} (x_i-x_j).$$
   Let us consider the case of $N=2$ in more detail. The group $G(1,1,3)$ consists of
   $6$ matrices:
   $$A_1=\left(\begin{array}{ccc}
      1 & 0 & 0\\
      0 & 1 & 0\\
      0 & 0 & 1
      \end{array}\right),\;\;
     A_2=\left(\begin{array}{ccc}
      0 & 1 & 0\\
      0 & 0 & 1\\
      1 & 0 & 0
      \end{array}\right),\;\;
     A_3=\left(\begin{array}{ccc}
      0 & 0 & 1\\
      1 & 0 & 0\\
      0 & 1 & 0
      \end{array}\right),$$
   $$A_4=\left(\begin{array}{ccc}
      0 & 1 & 0\\
      1 & 0 & 0\\
      0 & 0 & 1
      \end{array}\right),\;\;
     A_5=\left(\begin{array}{ccc}
      0 & 0 & 1\\
      0 & 1 & 0\\
      1 & 0 & 0
      \end{array}\right),\;
     A_6=\left(\begin{array}{ccc}
      1 & 0 & 0\\
      0 & 0 & 1\\
      0 & 1 & 0
      \end{array}\right).$$
   The reflection hyperplanes are eigenspaces of $1$ of matrices $A_4, A_5$ and $A_6$.
   These are thus three lines defined by linear factors of $(x_0-x_1)(x_0-x_2)(x_1-x_2)$.
   The resulting arrangement consists of $3$ lines passing through the point $(1:1:1)$.
   We refer to \cite[Example 1.3]{Sta07} for more details on these arrangements.
\end{example}
   In the next section we define hyperplane arrangements in general and we recall also
   some of their fundamental properties.

\section{Arrangements and their basic properties}
   Hyperplane arrangements are a classical subject of study in mathematics. It seems that the first
   non-trivial line arrangement studied by Greek geometers is that associated to the Theorem of Pappus.
   Nowadays hyperplane arrangements constitute an area of intensive study with far reaching connections
   to algebra, analysis, combinatorics, geometry and topology. In this section we establish the basic terminology
   and properties.
\begin{definition}[Hyperplane arrangement]
   A \emph{hyperplane arrangement} $\calh$ in the projective space $\P^N$ is a finite collection
   of mutually distinct hyperplanes.
\end{definition}
\begin{remark}
   Projective arrangements of hyperplanes correspond to \emph{central} affine arrangements, that is arrangements
   where all hyperplanes pass through the origin.
\end{remark}
   A hyperplane $H$ in $\P^N$ is defined by a linear polynomial $f_H$, which is determined uniquely up to a non-zero
   multiplicative scalar. To any arrangement $\calh$ one can thus associate its \emph{defining polynomial}
   $$Q(\calh)=\prod_{H\in\calh}f_H,$$
   which again is defined up to a scalar. It defines a unique principal ideal $I(\calh)=\langle Q(\calh)\rangle$ in $R$,
   which we call the \emph{arrangement ideal}.

   A fundamental combinatorial object associated to an arrangement is its intersection lattice.
\begin{definition}[Intersection lattice]
   Let $\calh$ be an arrangement. The set $L(\calh)$ of all non-empty intersections of hyperplanes in $\calh$
   is the \emph{intersection lattice} of $\calh$. This set has a natural structure of a poset defined by reversed
   inclusion relation.
\end{definition}
   We shall now introduce the arrangements which are in our focus.
\begin{example}[Fermat arrangements]\label{ex: Fermat arrangement}
   It is natural to extend the arrangements in Example \ref{ex: braid} by allowing
   powers. More precisely, for a positive integer $n$, let
   $$F_{N,n}=\prod_{0\leq i<j\leq N}(x_i^n-x_j^n).$$
   This polynomial splits over complex numbers into linear factors of the type
   \begin{equation}\label{eq: Fermat factors}
   x_i-\eps^kx_j,
   \end{equation}
   where $\eps$ is a primitive root of unity of degree $n$ and $k\in[n]$.
   The Fermat arrangement in $\P^N$ consists of zeroes of all linear factors of $F_{N,n}$.
   Following Orlik and Terao \cite[Example 6.29]{OrlTer} we denote this arrangement by $\cala^0_{N+1}(n)$.
   Of course, for $n=1$ we recover a braid arrangement.
\end{example}
   It is well known that Fermat arrangements are reflection arrangements. Indeed,
   the group $G(n,n,N+1)$ contains reflections in all hyperplanes defined in \eqref{eq: Fermat factors}.
   Thus the example generalizes readily as follows.
\begin{example}[Extended Fermat arrangements]
   For the groups $G=G(n,p,N+1)$ with $p<n$, we obtain reflection arrangement with
   $$Q(\calh(G))=x_0\cdot\ldots\cdot x_N\cdot\prod_{0\leq i<j\leq N}(x_i^n-x_j^n).$$
   Thus the reflection hyperplanes are all those of the corresponding Fermat arrangement
   with the addition of coordinate hyperplanes. We call the resulting arrangement $\calh(G)$
   the extended Fermat arrangement. In the literature it can be also encountered under the name
   of Ceva arrangement.

   Again, following \cite{OrlTer} we denote this arrangement by $\cala_{N+1}^{N+1}(n)$.
   In turn, following Hirzebruch \cite{Hir83}, we introduce intermediate Fermat arrangements $\cala_{N+1}^{k+1}(n)$
   as defined by linear factors of polynomials
   $$F_{N,n,k}=x_0\cdot\ldots\cdot x_k\cdot\prod_{0\leq i<j\leq N}(x_i^n-x_j^n),$$
   for $k=0,\ldots,N+1$.
\end{example}
   We introduce now briefly some useful properties of arrangements.

   We denote by $\Der(R)$ the $R$-module of $\C$-linear derivation of $R$. It is a free $R$-module with basis $D_0,\ldots,D_N$,
   where $D_i$ stands as usual for the partial derivation $\partial/\partial x_i$ for $i=0,\ldots,N$. We say that
   a derivation $\theta\in\Der(R)$ is homogeneous of polynomial degree $d$ if
   $$\theta=\sum_{i=1}^N f_iD_i,$$
   with $f_i$ a homogeneous polynomial in $R$ of degree $d$. In this way $\Der(R)$ becomes a $\Z$-graded $R-module$ with
   $$\Der(R)=\bigoplus_{d\in\Z}\Der(R)_d,$$
   where $\Der(R)_d$ consists of all homogeneous derivations of polynomial degree $d$.

   We want to distinguish these derivations which keep the arrangement ideal invariant.
\begin{definition}[Module of derivations]
   Let $\calh$ be an arrangement of hyperplanes in $V$.
   The \emph{module of $\calh$-derivations of $\calh$} is defined by
   $$\Der(\calh):=\left\{\theta\in\Der(R):\; \theta(Q(\calh))\in I(\calh)\right\}.$$
\end{definition}
   All this leads to the following important notion.
\begin{definition}[Free arrangements]
   We say that an arrangement $\calh$ is \emph{free} if the module $\Der(\calh)$ is a free $R$-module.
\end{definition}
   It is difficult to decide in general if an arrangement is free. However for some classes
   of arrangements it is known. In particular we have the following result, see \cite{Ter80}.
\begin{theorem}[Freeness of reflection arrangements]
   Any reflection arrangement is free.
\end{theorem}

\section{Unexpected curves}
\label{JSsec: unexpected curves}
   The concept of \emph{unexpected curves} has been introduced in the ground breaking
   article \cite{CHMN} of Cook II, Harbourne, Migliore and Nagel. Initially it has been defined only for
   curves in $\P^2$ with strong constrains on the relation between the degree and the multiplicity
   of the unexpected curve. More precisely we had (see \cite[Definition 2.1]{CHMN}).
\begin{definition}[Unexpected plane curves]\label{def: unexpected curve}
   We say that a finite set $Z$ of reduced points in $\P^2$ \emph{admits an unexpected curve}
   of degree $m+1$, if for a general point $P$, the fat point scheme $mP$ (i.e. defined
   by the ideal $I(P)^m$) fails to impose independent conditions on the linear system
   of curves of degree $m$ vanishing at all points of $Z$. In other words, $Z$ admits an
   unexpected curve of degree $m+1$ if
   $$h^0(\P^2, \calo_{\P^2}(d)\otimes I(Z+mP))\; >\; \max\left\{h^0(\P^2, \calo_{\P^2}(d)\otimes I(Z))-\binom{m+1}{2},\;0\right\}.$$
\end{definition}
   Note, that it follows immediately from the definition (taking a projection from the point $P$), that an unexpected curve is rational.
   Note also, that it is irrelevant if the points in $Z$ impose independent conditions on curves of any degree or not.
   In particular they can be arranged in a special position. In fact, in all example discovered so far the points
   in $Z$ exhibit a lot of symmetries.

   Research in \cite{CHMN} has been motivated by the article \cite{GIV} by Di Gennaro, Illardi and Valles, where
   the existence of unexpected curves has been first observed. Incidentally, in the example studied in \cite{GIV}
   the set $Z$ is dual to the $B_3$ arrangement of lines. It is the arrangement associated to the Weyl group
   of a $B_3$ root system.
\begin{definition}[$B_3$--arrangement of lines]
   The $B_3$ arrangement is the reflection arrangement defined by the the group $G(2,1,3)$.
\end{definition}
   Thus, according to Remark \ref{rmk: GnpN hyperplanes}, the lines in the $B_3$ arrangement
   are described by linear factors of the polynomial
   $$x_0x_1x_2(x_0^2-x_1^2)(x_1^2-x_2^2)(x_2^2-x_0^2).$$
   In the notation of Example \ref{ex: Fermat arrangement} this is the arrangement $\cala_3^3(2)$.

\begin{figure}[H]
\centering
\begin{tikzpicture}[line cap=round,line join=round,>=triangle 45,x=1.0cm,y=1.0cm,scale=0.45]
\clip(-7,-7) rectangle (7,7); \draw [line width=0.3pt,domain=-8.:2.]
plot(\x,{(-0.--1.*\x)/1.}); \draw [line width=0.3pt,domain=-3.:2.5]
plot(\x,{(-0.--2.*\x)/-2.}); \draw [line width=0.3pt,domain=-8.:8.]
plot(\x,{(-12.-4.*\x)/-2.}); \draw [line
width=0.3pt,domain=-1.:1.78] plot(\x,{(-6.--5.*\x)/-1.}); \draw
[line width=0.3pt,domain=-3.3:5.] plot(\x,{(-4.--1.*\x)/-3.}); \draw
[line width=0.3pt,domain=-8.:4.9] plot(\x,{(--18.-4.5*\x)/-7.5});
\draw [line width=0.3pt] (0.,-3.9) -- (0.,8.); \draw [line
width=1.2pt,] (0.005557527706670098,1.332101727856274)--
(0.007934814585178657,1.3331168080070637)-- (0.01,1.334)--
(0.012200674579804852,1.3349257486377006)--
(0.04435661922284132,1.3486142993799834)--
(0.07683655323255811,1.361033891769418)--
(0.11159778170207525,1.372822480568722)--
(0.14368380476106135,1.3824803982100569)--
(0.1663342151883824,1.388877842376568)--
(0.18579825381076956,1.3936697503622943)--
(0.20266873398728938,1.3977453422058248)--
(0.22499803717845537,1.4022655440686491); \draw [line width=1.2pt,]
(0.22499803717845537,1.4022655440686491)--
(0.22778920007735112,1.4029324590975905)--
(0.23102497299553998,1.4035746735699042)--
(0.23566867764149155,1.4044144924952526)--
(0.24134980566579264,1.405303712533841)--
(0.24824126096483617,1.4065634409218413)--
(0.2552809196036441,1.4076255648568219)--
(0.2621229737894328,1.4087370899050573)--
(0.2684710168426736,1.4097251121701555)--
(0.2744238509898768,1.4106884338786563)--
(0.281167102949156,1.4116270550304997)--
(0.2867988298602023,1.412516275069088)--
(0.29285046623391436,1.4131831900980294)--
(0.2975929731063744,1.4137760034570883)--
(0.30463263174518235,1.4145911218257945)--
(0.3100420536465821,1.4151839351848683)-- (0.318,1.416)--
(0.3256281048784691,1.4166906691391432)--
(0.33434740136793994,1.41711057860181)--
(0.3424738844983533,1.4178268947440062)--
(0.34941474091065206,1.41832090587657)--
(0.3568990095687531,1.418592611999472)--
(0.36710033945586773,1.4189137192356291)--
(0.37473281145373316,1.4190619225753938)--
(0.381772470092533,1.4192101259151586)--
(0.38851572205181206,1.419234826471786)--
(0.3949625673315623,1.4193089281416684)--
(0.40207632764025225,1.419284227585041)--
(0.407930359560945,1.4193583292549232)--
(0.4156369332286551,1.419284227585041)--
(0.4232200041132648,1.418963120348884)--
(0.4317169955930882,1.4186914142259819); \draw [line width=1.2pt,]
(0.4317169955930882,1.4186914142259819)--
(0.4394729703740897,1.4183950075464524)-- (0.446,1.418)--
(0.45278657039625475,1.4175304880644912)--
(0.45987563014831667,1.4168882735921773)--
(0.4673104976931621,1.4160978557800987)--
(0.475091173030791,1.41530743796802)--
(0.4833905600575952,1.4142453140330393)--
(0.4910477326120799,1.4133560939944507)--
(0.49910011407261007,1.41229397005947)--
(0.5060903715981623,1.411207145567862)--
(0.5131547307935967,1.410021518849744)--
(0.5207378016782064,1.4087617904617438)--
(0.5285925786857161,1.407329158177351)--
(0.5364226551365996,1.4056248197700565)--
(0.5434129126621519,1.4040933852591542)--
(0.5512676896696628,1.4024137474084872); \draw [line width=1.2pt,]
(0.5512676896696628,1.4024137474084872)--
(0.5586284555446259,1.4005859062180555)--
(0.5659151197497082,1.3985851611312312)--
(0.5737945973138466,1.3965103143745248)--
(0.5825138938033164,1.3942131626081713)-- (0.59,1.392)--
(0.5967661149773221,1.389865864641739)--
(0.6055595131366726,1.387074701742836)--
(0.6143529112960245,1.3841847366174236)--
(0.6233933150216503,1.380973664255854)--
(0.6320385098412379,1.3776390891111474)--
(0.6370774233932274,1.3755395417978136)--
(0.6444628898248179,1.372476672776009)--
(0.6522682657190741,1.3692408998578123)--
(0.659110319904862,1.366474437515537)--
(0.6659276735340225,1.3634115684937325)--
(0.6732884394089855,1.3599534905658883)--
(0.6796858835754803,1.3568906215440837)--
(0.6860833277419751,1.3536548486258868)--
(0.6925301730217247,1.3503202734811801)--
(0.6988041144050825,1.347183302789493)--
(0.703849221656292,1.3442631171942738)--
(0.7106096189491798,1.3405558025497777)--
(0.717089631100888,1.3369419496190085)--
(0.724223875248682,1.3327673264058784)--
(0.7299561936905776,1.3293092429979871); \draw [line width=1.2pt,]
(0.7299561936905776,1.3293092429979871)--
(0.7363738980331317,1.325383851021462)--
(0.743103141421444,1.3209911504763028)--
(0.7495519996685764,1.3170034506906263)--
(0.7549416251601414,1.3133895977598569)--
(0.7596770186556204,1.3100872838748436)--
(0.7653158753837894,1.3059749684708755)--
(0.770923578207383,1.301987268685199)--
(0.777154359122487,1.2974076447125862)--
(0.7825751385186274,1.2930149441674268)--
(0.7872482242049553,1.28952570685496)--
(0.7925443879827937,1.2850706985006493)--
(0.7979963212835096,1.2809272291920948)--
(0.8037286397254052,1.2761606817920284)--
(0.8100840362588112,1.2707087484912987)--
(0.8161590476510375,1.2656618159500517)--
(0.8203648247687327,1.26192334740098)--
(0.8259413736877508,1.2567829531460062)--
(0.8317983077479485,1.2516737127956084)--
(0.8354433145832842,1.248059859864839)--
(0.8402721697924899,1.2434490819876505)--
(0.845381410142875,1.2381529182097988)--
(0.8511448824893462,1.2325140614816157)--
(0.8565968157900621,1.2267817430397057)--
(0.8614256709992677,1.2214544253572783)--
(0.8672826050594654,1.2150990288238563)--
(0.8726410766464547,1.2089617096224636)--
(0.87775031699684,1.2033228528942805)-- (0.88,1.2)--
(0.8868986834066963,1.19058488931153); \draw [line width=1.2pt,]
(0.8868986834066963,1.19058488931153)--
(0.897380936434565,1.178047292552677)--
(0.9078631894624335,1.1642764895552482)--
(0.9175233049979202,1.1511222896771074)--
(0.9247170080562614,1.1406400366492138)--
(0.9312941079953162,1.1305688523675121)--
(0.9374601391881802,1.1207032024589065)--
(0.9458870484850941,1.1073434682076695)--
(0.9516948524706488,1.097767898134204)--
(0.9568961286433775,1.0887028168045774)--
(0.9618001890348074,1.0800835591468996)--
(0.9657686098888901,1.0725389255245832)--
(0.9705650291689956,1.0639161492906655)-- (0.975,1.055)--
(0.98,1.045); \draw [line width=1.2pt,] (0.98,1.045)--
(0.9835835373833466,1.0369252264131446)--
(0.9874786189553951,1.0285240700812552)--
(0.9905930364157587,1.0218906070079106)--
(0.9920453143191821,1.0186955956203716)--
(0.9933330958545789,1.0159785036293787)--
(0.9943000043479833,1.0139082517016138)--
(0.9954508034138442,1.0112032807014395)--
(0.9962376840029507,1.009310159677628)--
(0.996734312705552,1.008068587921122)--
(0.99719783282798,1.0069594504853099)--
(0.9980421016224025,1.0048405013542063)--
(0.9985828750985684,1.003510640050571)-- (0.999,1.0025); \draw [line
width=1.2pt,] (0.999,1.0025)--
(0.9993876171054071,1.0015231047009434)--
(0.9994769642680115,1.0012839375781342)-- (0.9996,1.001)--
(0.9997160179067515,1.0006937739074944)--
(0.9998243390868056,1.0004173681377009)--
(0.9999550715454915,1.0000849341713278); \draw [line width=1.2pt]
(-1.4579327674272156,-0.20799619816480194)--
(-1.196892165419736,-0.18051824005871248)--
(-0.8566503900477577,-0.1376849524227495)--
(-0.5883567394231182,-0.10088528009158544)--
(-0.4479690885763671,-0.07960953966472273)--
(-0.2605545190744047,-0.048970904872420624)--
(-0.09367063856469428,-0.018694277487243893)--
(3.036052656470094E-4,7.178715648877732E-5); \draw [line
width=1.2pt] (0.005557527706670098,1.332101727856274)--
(0.0026353136800025692,1.3310058975962709)--
(-0.0047371145375631986,1.3278246788787134)--
(-0.014722272305354273,1.3231402838765332)--
(-0.028035815995742373,1.3178395211109082)--
(-0.04795732732120207,1.3100755970419338)--
(-0.0725147911614706,1.2994124877428548)--
(-0.12583033765679041,1.2742087748541224)--
(-0.1910596993784429,1.239739269459747)-- (-0.26,1.2)--
(-0.31525163534964684,1.1650784281866293); \draw [line width=1.2pt]
(-0.31525163534964684,1.1650784281866293)--
(-0.37097411688511367,1.128658505614377)--
(-0.43325218448357666,1.0856829969791193)--
(-0.5071846273051438,1.032509910023631)--
(-0.6153517973445795,0.9480156896560026)--
(-0.7369943387357294,0.8467683049051415)--
(-0.8590010793526012,0.7371443379626611)--
(-0.9963041874497972,0.6085820112826112)--
(-1.1344788423768943,0.47137696652538863)--
(-1.225892163864254,0.37724839786500375)--
(-1.3182105677425775,0.28040458203172314)--
(-1.411434054011865,0.18084551902554685)--
(-1.496340100042257,0.08528652397526854); \draw [line width=1.2pt]
(-2.,2.)-- (-2.0061010188881174,1.9253915519147102)--
(-2.0126412894781573,1.8349178087523623)--
(-2.018091514969857,1.7291834342132328)--
(-2.0213616502648772,1.6299893302641528)--
(-2.022451695363217,1.5318852714134141)--
(-2.023541740461557,1.4250608517759433)--
(-2.022451695363217,1.3280468380235462)--
(-2.020271605166537,1.212502057599343)--
(-2.018091514969857,1.098047322273481)--
(-2.0104611992814774,0.9857726771443025)--
(-2.0006507933964173,0.86368762613005)--
(-1.9875702522163372,0.7307021241323821)--
(-1.971565998789021,0.6045255139060233)--
(-1.9535937676221757,0.47682808193088827)--
(-1.9318379088412576,0.32737479117480434)--
(-1.9006229810251576,0.1703542451905643)--
(-1.871299867016094,0.029414116566156074); \draw [line width=1.2pt]
(-2.,2.)-- (-1.9915032644961579,2.084383751349717)--
(-1.9833271569059092,2.1736920034894816)--
(-1.9738931866094687,2.262371324276149)--
(-1.9606856281944516,2.3604846153592707)--
(-1.9462202070732426,2.4579689750892895)--
(-1.9317547859520336,2.5497929526414413)--
(-1.9166604334777284,2.6365854793688177)--
(-1.9009371496503273,2.7246358688023884)--
(-1.8799591206362682,2.8330274741138224)--
(-1.856891311291082,2.9452908129272193)--
(-1.8245963782078216,3.079084107129487)--
(-1.8015285688626357,3.175968906379405)--
(-1.7707714897357207,3.286694391236454)--
(-1.7338629947834232,3.3974198760935033)--
(-1.693878791918434,3.520448192601336)--
(-1.6600460048788277,3.6188708458076015)--
(-1.6246753638828757,3.7126799371448236)--
(-1.5800775991488494,3.8234054220018727)--
(-1.5324041265021315,3.9479715924660574)--
(-1.483192799899068,4.055621369410411)--
(-1.4401328891213871,4.147892606791285)--
(-1.387845854605632,4.247853113953899)--
(-1.3386345280025684,4.341662205291121)--
(-1.2894232013995048,4.437009150584692)--
(-1.2371361668837497,4.524666826096523)--
(-1.1833112784116488,4.615400209521049)--
(-1.131024243895894,4.699982177120184)--
(-1.0618208158603357,4.807631954064537)-- (-1.,4.9)--
(-0.9326410835272935,4.989098720913587)--
(-0.8803540490115384,5.0644535647746345)--
(-0.8126884749323259,5.150573386330117)--
(-0.7496364627221507,5.22900393810386)--
(-0.6327595620398746,5.367410794175172)--
(-0.5758589656550823,5.428924952429088)--
(-0.5158826613575985,5.496590526508396)--
(-0.46051991892915195,5.555028976849616)--
(-0.40054361463166815,5.61808098905988)--
(-0.3359537484651472,5.68420870918284)--
(-0.27751529812400916,5.739571451611365)--
(-0.21600113987017963,5.801085609865281)--
(-0.14064629600923848,5.868751183944589)--
(-0.06007250391701196,5.94316180229321)--
(-0.03801637879127646,5.962953579003855); \draw [line width=1.2pt]
(-0.03801637879127646,5.962953579003855)--
(-0.016667409910643367,5.983278631602183)-- (0.,6.); \draw [line
width=1.2pt] (0.,6.)-- (0.007440940058896936,6.006056733359961)--
(0.02369806458035484,6.018249576751073)--
(0.015815822388132825,6.012214735072644)--
(0.030471866464295633,6.024653898532262)--
(0.03823094862226417,6.031550860450466)--
(0.04697531105426047,6.039309942608445)--
(0.10654089183360715,6.092910903332557); \draw [line width=1.2pt]
(0.10654089183360715,6.092910903332557)--
(0.33600770332278884,6.293800992945117)--
(0.6469746574028202,6.568265175603173)--
(0.8864238567355238,6.781756386512943)--
(1.2016280018748116,7.072591807826762)--
(1.4834575904699323,7.348064338032917)--
(1.7475845185077155,7.617749970021379)--
(1.902982648729627,7.782253459436214)--
(2.042598156350875,7.9352234938736395)--
(2.1967822386804277,8.108225318534998)--
(2.3855666859422033,8.32311179548282); \draw [line width=1.2pt]
(-2.421731464592546,-0.31500782845340164)--
(-3.0250476480032638,-0.3718223569361096)--
(-3.349702096475399,-0.4015823480460995)--
(-3.7365819809046936,-0.436753246630633)--
(-4.0720582443259,-0.46651323774062287)--
(-4.423767230170713,-0.4989786825878846)--
(-4.813352568337275,-0.534149581172418)--
(-5.140712470546678,-0.5639095722824079)--
(-5.484305095179688,-0.5963750171296696)--
(-5.82519226607543,-0.6261350082396595)--
(-6.1417303533357614,-0.6558949993496493)--
(-6.477206616756968,-0.6829495367223675)--
(-6.828915602601781,-0.7127095278323573)--
(-7.180624588446594,-0.7397640652050754)--
(-7.535039028028675,-0.7695240563150653)--
(-7.911097097508898,-0.801989501162327)--
(-8.24116245345557,-0.8263385847977732)--
(-8.563111448190437,-0.8533931221704913); \draw [line width=1.2pt]
(-2.0109660587517055,-0.5151610559327003)--
(-2.0558826671296444,-0.5702165047296889)--
(-2.094592199432454,-0.6218292144668616)--
(-2.1333017317352634,-0.6669903354868877)--
(-2.169246297445015,-0.7121514565069137)--
(-2.2079558297478243,-0.7609991996510236)--
(-2.2533393967777364,-0.8180001288433829)--
(-2.2945651640560945,-0.8701846443856862)--
(-2.336834621645297,-0.9213254696171436)--
(-2.3791040792344997,-0.9740318303148698)--
(-2.423460917445391,-1.0309129522559806)--
(-2.5236551872864634,-1.157199479868335)--
(-2.628425682633679,-1.2946819678289427)--
(-2.7796334558870357,-1.4923321285818356)--
(-2.945962006465728,-1.7083432332298059)--
(-3.1004099462888,-1.9157138936918574)--
(-3.262072970033177,-2.1330197218031617)--
(-3.3661036050054247,-2.274292327407467)--
(-3.484909712561634,-2.4384607669399228)-- (-3.6,-2.6)--
(-3.725762094243768,-2.763557479435115)--
(-3.859688979125313,-2.9514871404788474)--
(-3.9863246867435835,-3.1288872847188864)--
(-4.113533575947183,-3.304133097118505)--
(-4.233625883936595,-3.473151900955701)--
(-4.371509644961476,-3.6679683116943633)--
(-4.504945542727486,-3.857447286522378)--
(-4.681080927778624,-4.1029693384121995)--
(-4.831418705928333,-4.323583356052334)--
(-4.991541783247549,-4.549534809603111)--
(-5.181910330727059,-4.8199648957426255)--
(-5.393190755487162,-5.123365392875917)--
(-5.467151397253818,-5.231090675449246)--
(-5.538700278962864,-5.333188517888146)--
(-5.61587660080633,-5.44412948053829)--
(-5.707523482995434,-5.575972363687731)--
(-5.796758605126942,-5.707815246837177)--
(-5.891621167392868,-5.842873810063439)--
(-5.958346528986699,-5.940148132387116)-- (-6.,-6.); \draw [line
width=1.2pt] (-6.,-6.)-- (-6.1454624122092385,-6.211324631902503)--
(-6.271532861591601,-6.3963371095677966)--
(-6.415613375171445,-6.602634208557416)--
(-6.556419331624474,-6.808931307547034)--
(-6.7414318092895,-7.08235682763645)--
(-6.957617103723639,-7.4007595645619615)--
(-7.323035817220222,-7.940187189248189)--
(-7.876393949596476,-8.753616065345515)--
(-8.377242704232836,-9.491814669233868); \draw [line width=1.2pt]
(-1.7805615966263373,-0.31883845947713346)--
(-1.7728575551123,-0.34329779264779836)--
(-1.7637310132304043,-0.37253997786127246)--
(-1.7551632392188286,-0.39824329989604584)--
(-1.7443603937259724,-0.43046558041789945)--
(-1.733930060146663,-0.4611978132855631)--
(-1.7249897742215405,-0.4855973436229204)--
(-1.7147456965990056,-0.5131632252254312)--
(-1.7029510756722306,-0.5450676309058872)--
(-1.6877385828668638,-0.5859893043239564)--
(-1.6760044575901407,-0.6152314895374305)--
(-1.6633544518415526,-0.6462810319357901)--
(-1.647303328447496,-0.687185507681995)--
(-1.6307344268794375,-0.7270544270801947)--
(-1.616236638007386,-0.7638166774343789)--
(-1.5981144019173221,-0.8062744877025915)--
(-1.5774032749572489,-0.8539100797108301)--
(-1.5539103713977513,-0.9079445221829706)--
(-1.53483945596237,-0.9493324663193912)--
(-1.5141454838941901,-0.9927492312468128)--
(-1.496697625091607,-1.0304852979594317)--
(-1.4715402472832315,-1.0828288743672576)--
(-1.4498371756034836,-1.1242580763242478)--
(-1.4241965694704295,-1.1725227466924206)--
(-1.3970476923883723,-1.22267275574685)--
(-1.3698988153063152,-1.2705603583777714)--
(-1.342749938224258,-1.3180708932714413)--
(-1.322809517759893,-1.351123675369965)--
(-1.3069438637870714,-1.3772758522482809)--
(-1.288462992126422,-1.4069149860437056)--
(-1.2698077726199175,-1.435682380609853)--
(-1.2492347268090058,-1.4672393407096875)--
(-1.2235704651450579,-1.5057488372600978)--
(-1.2036101240199142,-1.534682175771816)--
(-1.1823679261252658,-1.5645311262617532)--
(-1.157463280317747,-1.598775014247141)--
(-1.1292624313886463,-1.6366813501454094)--
(-1.107653988702711,-1.6645159542832648)--
(-1.086594913203706,-1.6912518240472576)--
(-1.063704613748266,-1.7198189177676881)--
(-1.0360531320060944,-1.7527809489835697)--
(-1.0062041815162008,-1.787574204155889)--
(-0.9820320252912561,-1.8150425635024572)--
(-0.9532818091752233,-1.8461733707619008)--
(-0.9263628170156257,-1.8749235868779752)--
(-0.8932176634041477,-1.9082518628851444)--
(-0.8664817936401936,-1.9349877326491371)--
(-0.8375484551285174,-1.962456091995705)--
(-0.8018395879780309,-1.9948687560246607)--
(-0.767412577597049,-2.0250839513058856)--
(-0.7403104630418079,-2.048157373157003)--
(-0.712260074903584,-2.070395205919518)--
(-0.6896442378245455,-2.0887168967177527)--
(-0.6690323356765611,-2.1047483761662082)--
(-0.6438400108290245,-2.1233563433831653)-- (-0.62,-2.14)--
(-0.5931690847152294,-2.158568342886023)--
(-0.569746277555396,-2.1742575033846205)--
(-0.5469682226226843,-2.1889005386985274)--
(-0.5247818054804327,-2.202951936221974)--
(-0.49386873092889544,-2.221144798278646)--
(-0.4659138453296584,-2.2371190186210903)--
(-0.44210042426364166,-2.250282959458845)--
(-0.41991400712139004,-2.2616719869252173)--
(-0.3925507593126129,-2.275575475001064)--
(-0.366814515427601,-2.287851959153128)-- (-0.34,-2.3)--
(-0.31741275992418744,-2.309742557400181)--
(-0.2953742522295508,-2.3190608525999403)--
(-0.2768855712776745,-2.326160506085471)--
(-0.25558661082111295,-2.333999706809078)--
(-0.2316252803074812,-2.342578454770761)-- (-0.21,-2.35)--
(-0.1893231782895881,-2.356481942846592)--
(-0.1615162021379661,-2.3649127813607937)--
(-0.13193431261496394,-2.3730478009796316)--
(-0.10338778922526687,-2.3799995450175473)--
(-0.07365799025464971,-2.386655470160233)--
(-0.04940084084578795,-2.391684391379151)--
(-0.029285155970146484,-2.3952342181219164)--
(-0.017181955492201847,-2.3970681960307965)--
(-0.00751229028756571,-2.3986985465594874)-- (0.,-2.4)--
(0.01246138989065485,-2.4016990969121728)--
(0.023423829207000698,-2.4030568302219972)--
(0.03343082656458246,-2.404263704275175)--
(0.04816474729710234,-2.4059231560982943)--
(0.06475926552826808,-2.4076328943402956)--
(0.0783365986264946,-2.4087391955557087)--
(0.09095848976595702,-2.409644351095592)--
(0.10398267225647802,-2.4104489337977104)--
(0.11861602015123444,-2.411152943662064)--
(0.13093619277740307,-2.4116558078508876)--
(0.14320607898468937,-2.412058099201947)--
(0.15633083431297515,-2.4122592448774767)--
(0.172573347611963,-2.412309531296359)--
(0.19052559915295136,-2.4121586720397117)--
(0.20968472474711544,-2.4117060942697703)--
(0.22617867014051654,-2.4110020844054167)-- (0.245,-2.41)--
(0.26383047326146847,-2.408524885856892)--
(0.28197433832250585,-2.407129203929118)--
(0.3023668020449538,-2.404880605267704)--
(0.3250078644288124,-2.402011703527279)--
(0.35960070644029346,-2.396862706207764)--
(0.38908959088648465,-2.391446380493149)--
(0.4251984289838615,-2.384024008217566)--
(0.4619090810495281,-2.3751974033493046)--
(0.5018294076127393,-2.3641641472639776)--
(0.5493727111076189,-2.3489181934006176)--
(0.5977184332268834,-2.3314655883201856)--
(0.6462647600022458,-2.311405122710501)--
(0.6982213659312491,-2.2873325639788793)--
(0.7455640647700321,-2.2630594005911604)--
(0.795314019481956,-2.234974748737602)--
(0.8524863464694681,-2.199267119952363)--
(0.908053836208208,-2.161352839950059)--
(0.9700406749420382,-2.114411350423402)--
(1.0234015134637162,-2.070679535394289)--
(1.0693399797098224,-2.0305586041749195)--
(1.1168832832047029,-1.9850213472409217)--
(1.1670344472288365,-1.935070787872807)--
(1.1955203083945443,-1.9045788801460863)--
(1.219793471782225,-1.8768954376047216)--
(1.2518939974923946,-1.8399614077289006)--
(1.2810538404037113,-1.8048615968170758)--
(1.308053694951227,-1.7719217742690554)--
(1.3369435393170686,-1.734661974993426)--
(1.3596234171369819,-1.704152139354686)--
(1.3836532876842706,-1.6720223124430922)--
(1.40876315241346,-1.6366525029857917)--
(1.4333330200516994,-1.6015526920739669)--
(1.4519629196894852,-1.5742828389809336)--
(1.4714028149636964,-1.5448529975240959)--
(1.490032714601482,-1.5154231560672582)--
(1.5092880260065173,-1.4859001832953334)-- (1.5,-1.5)--
(1.520389210688688,-1.4680023549301733)--
(1.5360214911594996,-1.4419485541454464)--
(1.5625284015230498,-1.3989031441532889)--
(1.5867697639922793,-1.3592560560026175)--
(1.6108298694489156,-1.3198887199809326)--
(1.6427049500444193,-1.263638577753485)--
(1.675517533010379,-1.2087946890817236)--
(1.7036426041240589,-1.1591070634474783)--
(1.7284864169411427,-1.115981954406435)--
(1.7547364833139105,-1.069106835883562)--
(1.7833303056128185,-1.020825463805003)--
(1.813799132652641,-0.9702003358003002)--
(1.833486682432217,-0.9378565040195178)--
(1.8522367298413367,-0.906450174609193)--
(1.8705180260652288,-0.8787938546806978)--
(1.8892680734743486,-0.8511375347522028)--
(1.9070806185130125,-0.8239499660089364)--
(1.9248931635516764,-0.7976998996361275)--
(1.9435964809711592,-0.7707398139943259)--
(1.9633273325837743,-0.7437608944423406)--
(1.9832595194169262,-0.7169833101108924)--
(2.0019836949268517,-0.6940310949696444)--
(2.0186945182314138,-0.6734949024748499)--
(2.036009347197586,-0.6539653860827412)--
(2.0513108239583895,-0.6370532275576163)--
(2.072249686894226,-0.6145036828574496)--
(2.096208578138116,-0.589940785951911)--
(2.1298315599677817,-0.5581298211070334)--
(2.1630518713563682,-0.528936214129139)--
(2.1960708475244184,-0.5017559593566168)--
(2.2228484318558244,-0.48162243730289656)--
(2.2522433740542045,-0.46048223914648984)--
(2.2834503332374223,-0.43954337621062084)--
(2.309019906245607,-0.4234365585676447)--
(2.3353948201359396,-0.40813508180681735)--
(2.3625750749084196,-0.39283360504599)--
(2.391164676224658,-0.3771294578440883)--
(2.418344930997138,-0.36303599240648415)--
(2.450357231062503,-0.3477345156456568)--
(2.484584218553774,-0.3322317036642923)--
(2.524851262661152,-0.3147155394775556)--
(2.56572231243014,-0.2986087218345794)--
(2.6051840156553703,-0.2839112507353637)--
(2.653303133363687,-0.26659642176916437)--
(2.7133010290836794,-0.24686557015651858)--
(2.7765202883322626,-0.22753738898494724)--
(2.826048752584337,-0.21324258832680593)--
(2.866315796691715,-0.202370486417797)--
(2.9087975282249987,-0.1919010549498625)--
(2.977855508869151,-0.17438489076312597)--
(3.025571956136394,-0.16351278885411705)--
(3.0805364713429646,-0.15143267562188495)--
(3.159057207352351,-0.1351245227583716)--
(3.232343227628072,-0.12042705165915586)--
(3.3100586227553186,-0.10613225100101453)--
(3.3797206090610885,-0.09344813210717082)--
(3.4487785897052476,-0.08177068931601311)--
(3.512400519394947,-0.07150259306861582)--
(3.588505232757899,-0.05922114461584651)--
(3.6712540083985687,-0.046537025722002794)--
(3.74152000036595,-0.0360675942540683)--
(3.8115846571127645,-0.02620216844774541)--
(3.88547468304981,-0.015531401759273714)--
(3.956747351119875,-0.005867311173488026)-- (4.,0.)--
(4.056609620506188,0.007622148602504496)--
(4.1365397030593405,0.017487574408827387)--
(4.247676744795713,0.03198371028750592)--
(4.375725945057187,0.04748652226887046)--
(4.51223122458121,0.0633920046913094)--
(4.629408322933691,0.07688146446730193)--
(4.905640245510354,0.1070817475478822)--
(5.273133169494088,0.14559636797820494)--
(5.643348934558216,0.18261794448467142)--
(6.027535105851179,0.22033804130258067)--
(6.497639275451477,0.2657418615463603)--
(6.927229266988155,0.3048589989871551)--
(7.514684848910693,0.3579465426568051)--
(8.098716413217973,0.4088154119643845); \draw [line width=1.2pt]
(3.036052656470094E-4,7.178715648877732E-5)--
(0.0013055303442626436,2.8052154786733233E-4)--
(0.002630458171969699,5.856850042979952E-4)--
(0.005358076966663793,0.001155865204319609)--
(0.010352705921614642,0.002299448619662259)--
(0.022030131842021718,0.0047274282664830376)--
(0.03959551073952203,0.008608414235823985)--
(0.061281302413794314,0.013283532960414035)--
(0.09473938099695726,0.020662334802839296)--
(0.12684561802120423,0.027984809913642994)--
(0.15562857787978374,0.03480034443985259)--
(0.1921282999705071,0.043587314572817025)--
(0.2177569628582834,0.05000856197767565)--
(0.255721179971165,0.05975308654820672)--
(0.2843914863665016,0.06741352204873981)--
(0.31244219871400214,0.07524293774413761)--
(0.34415414861689897,0.0842552148035883)--
(0.36730443531333073,0.09123972952466254)--
(0.3979461773154194,0.10058996697384265)--
(0.4312352757037172,0.11134837271356124)--
(0.4605251761468902,0.12120555074733538)--
(0.4947155022411325,0.13325947131435062)--
(0.5183164027905355,0.1419901147156934)--
(0.5436634320202044,0.15184729274946754)--
(0.5692357681763596,0.16215508463621423)--
(0.5940195300898136,0.17257552998620404)--
(0.6183526781502958,0.1833902624575448)--
(0.6389119351921383,0.19290948010158968)--
(0.659020578381009,0.20276665813536385)--
(0.679974122544202,0.21341241041183992)--
(0.7034060428987403,0.22580429136858432)--
(0.723176725697882,0.23690065749803293)--
(0.7442429233243172,0.2494615186496438)--
(0.7645768734396738,0.26224768672773946)--
(0.7831646948747645,0.2745832409528596)--
(0.8008512886039112,0.28703144864122576)--
(0.8164537932630638,0.2985221018920256)--
(0.8342530404554536,0.31243480460255246)--
(0.8505877926256846,0.32600954692334955)--
(0.8670351982591586,0.3407671506081999)--
(0.8816238217491227,0.35451087312386215)--
(0.8961561185074662,0.3689868431506047)--
(0.9097308608282444,0.3836317933722119)--
(0.9225733556379433,0.3985020505203054)--
(0.9337823752305338,0.41218944630434606)--
(0.9461742561872596,0.42852419847460027)--
(0.9564257213423699,0.4430564952329645)--
(0.9642551370377558,0.45522306926322276)--
(0.9726663621557098,0.46882481975162216)--
(0.9789940136402842,0.4793709055592628)--
(0.9878869292402266,0.49601661937456587)--
(0.9965518213632428,0.5134034094898653)--
(1.00470366066319,0.5314172641667001)-- (1.01,0.545)--
(1.0163328579862076,0.5614593572511682)--
(1.021520392086214,0.57662291846648)--
(1.0253967911938815,0.5892212155664182)--
(1.0288171433477313,0.6014204715817417)--
(1.0318384544169428,0.6132206865124531)--
(1.0353728183092352,0.6286122712046854)--
(1.0391922115476804,0.6481652843507464)--
(1.0419284932707398,0.6656090803352763)--
(1.0442359008877171,0.6842627407300443)--
(1.04583193080892,0.6997933395786954)--
(1.0469982603667218,0.7160605676217502)--
(1.0475507322625228,0.7279080204983891)--
(1.0478576610935233,0.7462623645922495)--
(1.047734889561123,0.7647394802185087)--
(1.047182417665322,0.7795334498727573)--
(1.046200245406122,0.7964145355778112)--
(1.0447269870173195,0.8133570070490669)--
(1.0432694280382249,0.8269089584957893)--
(1.0413821226564024,0.8414266922021371)--
(1.0390592852633902,0.8558718372399531)--
(1.0360831498535945,0.8730027630134419)--
(1.0328892484382026,0.8887545040848293)--
(1.029695347022811,0.9033448264597087)--
(1.026574034275947,0.9162656094583579)--
(1.0238882535402767,0.9262828457157378)--
(1.0209121181304799,0.9370259686584351)--
(1.0171375073668352,0.9494386309773625)--
(1.014117818755902,0.95918003029433)--
(1.0110037648758952,0.9685221919343648)--
(1.0078098634605035,0.9777046585036299)-- (1.005,0.985)--
(1.0025399261251071,0.9921570624082992)--
(1.0000956830303942,0.9997304687014905)--
(1.0000581673361866,0.9998272762384144)--
(1.0000184252260595,0.9999291790848941)--
(0.9999847972867213,1.0000117203905425)--
(0.9999550715454915,1.0000849341713278); \draw [line width=1.2pt]
(-1.7805615966263373,-0.31883845947713346)--
(-1.787266990469414,-0.297649414932973)--
(-1.792094874036429,-0.28102003820211285)--
(-1.7970568654803059,-0.26519530873242336)--
(-1.8015368323981797,-0.2511411586067161)--
(-1.8060334998111258,-0.23393377682782404)--
(-1.8083471699002662,-0.22511735359381138)--
(-1.811825653588652,-0.2124200627969742)--
(-1.8152451457855412,-0.1988649283971524)--
(-1.8188080757435487,-0.1835652879892094)--
(-1.8221687853342348,-0.16897890946545055)--
(-1.8254785874928212,-0.15449852502160932)--
(-1.8291019844924927,-0.1400193068223624)--
(-1.833100537387191,-0.1233586697610905)--
(-1.837432303023114,-0.10536518173491682)--
(-1.8414308559178123,-0.08803811919119403)--
(-1.8450804697244878,-0.07286501419684081)--
(-1.8487857481977452,-0.05804390030378533)--
(-1.8524910266710024,-0.044281437403090956)--
(-1.8561963051442598,-0.030518974502396573)--
(-1.859901583617517,-0.01622718610552164)--
(-1.8633467225857363,-0.002389405080179365)--
(-1.8666138055378099,0.010025510137722708)--
(-1.869322799897498,0.020402262952604503)--
(-1.871299867016094,0.029414116566156074); \draw [line width=1.2pt]
(-2.421731464592546,-0.31500782845340164)--
(-2.3166344405164723,-0.30331876168655963)--
(-2.170305625126053,-0.28714557682759007)--
(-2.010884231516491,-0.26789178532881674)--
(-1.877647994345215,-0.25248875212979804)--
(-1.7459520604938379,-0.23785587059073032)--
(-1.6096352166827632,-0.2232229890516626)--
(-1.4579327674272156,-0.20799619816480194); \draw [line width=1.2pt]
(-2.0109660587517055,-0.5151610559327003)--
(-1.963660956215121,-0.4527839017255142)--
(-1.919502145095105,-0.3952941287578435)--
(-1.8669432832706956,-0.3342726887883281)--
(-1.82362547331279,-0.2817383235201345)--
(-1.7747777301687684,-0.22183071400377338)--
(-1.7213217093696507,-0.1610014489563913)--
(-1.67155231069461,-0.10662377262615581)--
(-1.6236262230816079,-0.055011062888983144)--
(-1.5626074124038152,0.012619097182656953)--
(-1.496340100042257,0.08528652397526854);
\begin{scriptsize}
\draw [fill=black] (0.,0.) circle (3.pt);
\draw [fill=black] (-2.,2.) circle (3.pt);
\draw [fill=black] (1.,1.) circle (3.pt);
\draw [fill=black] (0.,6.) circle (3.pt);
\draw [fill=black] (-6.,-6.) circle (3.pt);
\draw [fill=black] (1.5,-1.5) circle (3.pt);
\draw [fill=black] (4.,0.) circle (3.pt);
\draw [fill=black] (0.,1.33) circle (3.pt);
\draw [fill=black] (0.,-2.4) circle (3.pt);
\draw [fill=white] (-1.8035052604104964,-0.24410849141067362) circle
(6pt); \draw [fill=black] (-1.8035052604104964,-0.24410849141067362)
circle (3pt); \draw[color=black] (-2.3,0.1) node {$P$};
\end{scriptsize}
\end{tikzpicture}
\caption{A visualization of an unexpected quartic admitted for
$B_3$} \label{fig: B3}
\end{figure}
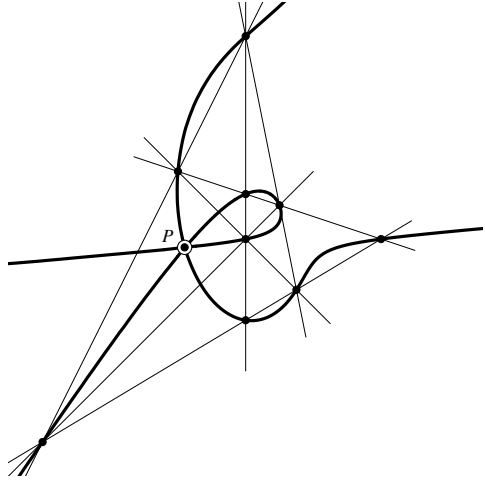

   Dually we obtain a set $Z$ of $9$ points with the following coordinates
   $$\begin{array}{lll}
   P_1=(1:0:0),& P_2=(0:1:0),& P_3=(0:0:1),\\
   P_4=(1:1:0),& P_5=(1:-1:0),& P_6=(1:0:1),\\
   P_7=(1:0:-1),& P_8=(0:1:1),& P_9=(0:1:-1).
   \end{array}$$
   Figure \ref{fig: B3} shows an unexpected curve admitted by a $B_3$ arrangement. The coordinate
   system in this Figure has been so chosen that the set $Z$ is completely contained in the affine part
   of the plane.
   It has been realized in \cite{BMSS} that one can actually write explicitly the equation of an unexpected quartic $Q_P$ in this case.
   If $P=(a:b:c)$ is general, then
   \begin{equation}\numberthis\label{eq: B3 quartic}
\begin{array}{rcl}
      Q_P(x:y:z) & = & 3a(b^2-c^2)\cdot x^2yz +3b(c^2-a^2)\cdot xy^2z
                      +3c(a^2-b^2)\cdot xyz^2\\
                     && +a^3\cdot y^3z -a^3\cdot yz^3 +b^3\cdot xz^3
                      -b^3\cdot x^3z +c^3\cdot x^3y -c^3\cdot xy^3.
\end{array}
   \end{equation}
   It is natural to wonder if the set of points dual to the Fermat arrangement $\cala_3^0(2)$ also admits an unexpected curve.
   It turns out that this does not happen for degree $2$ but allowing high enough degree we obtain in this
   way additional examples of unexpected curves. More precisely we have the following result, see \cite[Theorem 6.12]{CHMN}.
\begin{theorem}[Cook II, Harbourne, Migliore, Nagel]\label{thm: Fermat m>=5}
   For $m\geq 5$ let $Z$ be the set of points dual to lines in the Fermat arrangement $\cala_3^0(m)$, i.e., given
   by linear factors of
   $$F_{2,m}=(x_0^m-x_1^m)(x_1^m-x_2^m)(x_2^m-x_0^m).$$
   Then $Z$ admits an unexpected curve of degree $m+2$ with a point of multiplicity $m+1$.
   Moreover the unexpected curve is unique and irreducible.
\end{theorem}
   It is natural to wonder if for lower values of $m$ intermediate Fermat arrangements $\cala_3^k(m)$
   for $k=1,2$ might lead to unexpected curves. We show that this is indeed the case.
\subsection{Fermat-type arrangement for $m=3$}
\begin{theorem}\label{thm: Fermat m=3}
   Let $Z$ be the set of points dual to lines in the Fermat arrangement $\cala_3^2(3)$.
   Then $Z$ admits a unique and irreducible unexpected quintic with a point of multiplicity $4$
   at a general point $R=(a:b:c)$.
\end{theorem}
\proof
   The lines in $\cala_3^2(3)$ are given by linear factors of
   $$x_0 x_1 (x_0^3-x_1^3) (x_1^3-x_2^3) (x_2^3-x_0^3),$$
   so that the $11$ points in $Z$ are the coordinate points $(1:0:0)$ and $(0:1:0)$ together
   with the following points
      $$\begin{array}{lll}
      P_1=(1:-1:0),& P_2=(1:-\eps:0),& P_3=(1:-\eps^2:0),\\
      P_4=(1:0:-1),& P_5=(1:0:-\eps),& P_6=(1:0:-\eps^2),\\
      P_7=(0:1:-1),& P_8=(0:1:-\eps),& P_9=(0:1:-\eps^2).\\
   \end{array}$$
   The points $P_1,\ldots,P_9$ form a complete intersection given by the coordinate axes
   and the lines in the corresponding Fermat arrangement $\cala_3^0(3)$, i.e., their
   ideal is generated by $x_0x_1x_2$ and $x_0^3+x_1^3+x_2^3$. Intersecting with the ideals
   of the remaining $2$ points we obtain
   $$I(Z)=\langle x_0x_1x_2,\;
                  x_0^3x_2+x_1^3x_2+x_2^4,\;
                  x_1^4x_2+x_1x_2^4,\;  \rangle.$$
   Since the regularity of $I(Z)$ is $5$, the linear system of quintics vanishing in all points of $Z$
   has (projective) dimension $9$. Thus it is not expected that for a general point $R$,
   it contains a member which vanishes at $R$ to order $4$. However there exists such an unexpected
   curve given by the equation
   \begin{align*}
      Q_R(x_0:x_1:x_2) & ={} a^4\cdot x_1x_2\cdot (x_1^3+x_2^3)+b^4\cdot x_0x_2\cdot (x_0^3+x_2^3)+c^4\cdot x_0x_1\cdot (x_0^3+x_1^3)\\
                       & -4a(b^3+c^3)\cdot x_0^3x_1x_2
                        -4b(a^3+c^3)\cdot x_0x_1^3x_2\\
                       & -4c(a^3+b^3)\cdot x_0x_1x_2^3\\
                       & +6a^2b^2\cdot x_0^2x_1^2x_2+
                        6a^2c^2\cdot x_0^2x_1x_2^2+
                        6b^2c^2\cdot x_0x_1^2x_2^2.
   \end{align*}
   The latter claim can be easily verified by a direct computation.
\endproof
   Passing to the next degree of the Fermat-type arrangement we can drop another coordinate line.
\subsection{Fermat-type arrangement for $m=4$}
\begin{theorem}\label{thm: Fermat m=4}
   Let $Z$ be the set of points dual to lines in the Fermat arrangement $\cala_3^1(4)$.
   Then $Z$ admits a unique and irreducible unexpected sextic with a point of multiplicity $5$
   at a general point $R=(a:b:c)$.
\end{theorem}
\proof
   The lines in $\cala_3^1(4)$ are given by linear factors of
   $$x_0 (x_0^4-x_1^4) (x_1^4-x_2^4) (x_2^4-x_0^4),$$
   so that the $13$ points in $Z$ are the coordinate point $(1:0:0)$ and the points with coordinates
      $$\begin{array}{llll}
      P_1=(1:1:0),& P_2=(1:i:0),& P_3=(1:-1:0),& P_4=(1:-i:0),\\
      P_5=(1:0:1),& P_6=(1:0:i),& P_7=(1:0:-1),& P_8=(1:0:-i),\\
      P_9=(0:1:1),& P_{10}=(0:1:i),& P_{11}=(0:1:-1),& P_{12}=(0:1:-i).\\
   \end{array}$$
   The ideal of $Z$ is generated by
   $$x_0x_1x_2, x_0^4x_2+x_1^4x_2-x_2^5, x_0^4x_1-x_1^5+x_1x_2^4,$$
   so $Z$ is an almost complete intersection ideal (that means that the number of generators
   is one higher than the hight of an ideal). Such ideals have an easy minimal free resolution
   and either writing it explicitly down or using a symbolic algebra program (we used Singular \cite{Singular}) we get $\reg(I(Z))=6$, so that $Z$ imposes independent conditions
   on curves of degree $6$. Thus $\dim(I)_{[6]}=15$ and we do not expect that for a general point $R$
   there exists an element vanishing at $R$ to order $5$.

   However, such an element can be written down explicitly as follows:
   \begin{align*}
      S_R(x_0:x_1:x_2) & ={}  a^5\cdot x_1x_2\cdot (x_1^4-x_2^4)+ b^5\cdot x_0x_2\cdot (x_2^4-x_0^4)+ c^5\cdot x_0x_1\cdot (x_0^4-x_1^4)\\
      &+10a^3x_0^2x_1x_2\cdot (b^2x_1^2-c^2x_2^2)\\
      &+10b^3x_0x_1^2x_2\cdot (c^2x_2^2-a^2x_0^2)\\
      &+10c^3x_0x_1x_2^2\cdot (a^2x_0^2-b^2x_1^2)\\
      &+ 5a(b^4-c^4)\cdot x_0^4x_1x_2+
      5b(c^4-a^4)\cdot x_0x_1^4x_2+
      5c(a^4-b^4)\cdot x_0x_1x_2^4.\\
            \numberthis\label{eq:sextic}
   \end{align*}
   Vanishing order in $R$ can be checked by a direct (but tedious) computation.
   \endproof
   Theorems \ref{thm: Fermat m=3} and \ref{thm: Fermat m=4} fill thus a gap between the $B_3$ example
   and a general Theorem \ref{thm: Fermat m>=5}, showing that all these examples belong in fact
   to the same family. We summarize this section by the following statement.
\begin{theorem}\label{thm: Fermat m}
   Let $Z$ be the set of points dual to lines in the Fermat arrangement $\cala_3^k(m)$.
   Then for $m\geq 2$ and $0\leq k\leq 3$ such that $k+m\geq 5$, the set $Z$ admits
   a unique and irreducible unexpected curve $C_R$ of degree $m+2$ with a point of multiplicity $m+1$
   at a general point $R=(a:b:c)$.\\
   Moreover, the curve $C_R$ does not depend on $k$. Thus, for example an unexpected curve
   for $\cala_3^0(m)$ automatically passes through all three coordinate points.
\end{theorem}
\proof
   We provide general formulas for curves $C_R$ depending on the parity of $m$.
   We omit lengthy and not very instructive computational arguments showing that
   these formulas indeed define curves satisfying conditions required in our statement.

   For an even $m\geq 2$ we have the following formula.
   \begin{align*}
      C_R&(x_0:x_1:x_2) = \\
      & \sum_{k=1}^{\frac{m}{2}+1} \binom{m+1}{2k-1}a^{2k-1}\cdot  \Big[ \left(b^{m-(2k-2)}x_1^{2k-2}-c^{m-(2k-2)}x_2^{2k-2}\right)\cdot x_0^{m-(2k-2)}x_1x_2\\
      & +\left(c^{m-(2k-2)}x_2^{2k-2}-a^{m-(2k-2)}x_0^{2k-2}\right)\cdot x_0x_1^{m-(2k-2)}x_2\\
      & +\left(a^{m-(2k-2)}x_0^{2k-2}-b^{m-(2k-2)}x_1^{2k-2}\right)\cdot x_0x_1x_2^{m-(2k-2)} \Big].
   \end{align*}
   For an odd $m\geq 3$ we have in turn
    \begin{align*}
      C_R&(x_0:x_1:x_2) = \\
      & a^{m+1}x_1x_2(x_1^m+x_2^m)+b^{m+1}x_0x_2(x_0^m+x_2^m)+c^{m+1}x_0x_1(x_0^m+x_1^m)\\
      & -(m+1)\left[a(b^m+c^m)x_0^mx_1x_2+b(a^m+c^m)x_0x_1^mx_2+c(a^m+b^m)x_0x_1x_2^m\right]\\
      & +\sum_{k=2}^{\frac{m-1}{2}}(-1)^k\binom{m+1}{k}\Big[a^{m+1-k}x_0^kx_1x_2(b^kx_1^{m-k}+c^kx_2^{m-k})\\
      & +b^{m+1-k}x_0x_1^kx_2(a^kx_0^{m-k}+c^kx_2^{m-k})+
      c^{m+1-k}x_0x_1x_2^k(a^kx_0^{m-k}+b^kx_1^{m-k})\Big]\\
      & +\binom{m+1}{\frac{m+1}{2}}\cdot (-1)^{\frac{m+1}{2}}\Big[ a^{\frac{m+1}{2}}b^{\frac{m+1}{2}}x_0^{\frac{m+1}{2}}x_1^{\frac{m+1}{2}}x_2+b^{\frac{m+1}{2}}c^{\frac{m+1}{2}}x_0x_1^{\frac{m+1}{2}}x_2^{\frac{m+1}{2}}\\
      & +a^{\frac{m+1}{2}}c^{\frac{m+1}{2}}x_0^{\frac{m+1}{2}}x_1x_2^{\frac{m+1}{2}}\Big].
    \end{align*}
\endproof

\section{Unexpected hypersurfaces}
\label{JSsec:3}
   It has been quickly realized that Definition \ref{def: unexpected curve} is too restrictive.
   First, exactly the same phenomena can be studied in projective spaces of arbitrary dimension.
   More importantly,
   there is no need to couple the degree of the unexpected hypersurface and its multiplicity in a general point.
   In the subsequent to \cite{CHMN} article \cite{HMNT} Harbourne, Migliore, Nagel and Teitler
   generalize Definition \ref{def: unexpected curve} in the following way.
\begin{definition}[Unexpected hypersurface]\label{def: unexpected hypersurface}
   Let $Q_1,\ldots,Q_s$ be mutually distinct points in $\P^N$ and let $n_1,\ldots,n_s$
   be positive integers. Let $Z=n_1Q_1+\ldots+n_sQ_s$ be a scheme of fat points, i.e.,
   $$I(Z)=I(Q_1)^{n_1}\cap\ldots\cap I(Q_s)^{n_s}.$$
   Let $m$ be a positive integer and let $R$ be a general point in $\P^N$.
   We say that $Z$ admits an \emph{unexpected hypersurface} with respect to $X=mR$ of degree $d$, if
   $$h^0(\P^N, \calo_{\P^N}(d)\otimes I(Z+X))\;>\;\max\left\{h^0(\P^N, \calo_{\P^N}(d)\otimes I(Z))-\binom{N+m-1}{N},\; 0\right\}.$$
\end{definition}
\begin{remark}
   Note that in fact one can pose the same definition replacing a fat points scheme $Z$
   by an arbitrary subscheme of $\P^N$.
\end{remark}
   Interestingly, the first example of an unexpected surface in $\P^3$ has been announced
   by Bauer, Malara, Szemberg and the author in \cite[Theorem 1]{BMSS} and it is related
   to a Fermat-type arrangement of planes in $\P^3$.

   Let us begin with a general definition of a Fermat-derived configuration of flats.
\begin{definition}[Fermat-derived configurations]
   Let $\calh=\cala_{N+1}^{k+1}(n)$ be a Fermat-type arrangement of hyperplanes in $\P^N$.
   Taking intersections of hyperplanes in the arrangement we obtain a number of related
   objects. For an integer $0\leq t\leq N-1$ we denote by $\calh(t)$ the set theoretical union
   of all $t$--dimensional flats in the intersection lattice. We call these sets \emph{Fermat-derived
   configurations} of flats.
\end{definition}
   In particular $\calh(N-1)$ is the union of all arrangement hyperplanes, while
   $\calh(0)$ is the union of all points in $L(\calh)$. These configurations have been introduced
   in \cite{MalSzp18} and investigated further in \cite{MalSzp18Moieciu}
   in the context of the containment problem between symbolic and ordinary powers
   of homogeneous ideals, see \cite{SzeSzp17} for an introduction to this circle
   of ideas.

   We are now in a position to recall the main result from \cite{BMSS}.
\begin{theorem}[Unexpected quartic surface]\label{thm: BMSS}
   Let $Z$ be the subset of points $(\cala_{4}^{4}(3))(0)$ derived from the Fermat-type
   arrangement $(\cala_{4}^{4}(3))$ and defined by the following binomial ideal:
   $$x_0(x_1^3-x_2^3),\; x_0(x_2^3-x_3^3),\; x_1(x_0^3-x_2^3),\; x_1(x_2^3-x_3^3)\,,$$
   $$x_2(x_0^3-x_1^3),\; x_2(x_1^3-x_3^3),\; x_3(x_0^3-x_1^3),\; x_3(x_1^3-x_2^3)\,.$$
   Then $Z$ admits a unique and irreducible unexpected
   quartic surface $Q_R$, which vanishes at a general point $R=(a:b:c:d)$ to order $3$.
\end{theorem}
\proof
   Also in this case the statement is effective in the sense that we can write down
   explicitly the equation for $Q_R$:
   \begin{align*}
      Q_R(x_0:x_1:x_2:x_3)&=b^2(c^3-d^3)\cdot x_0^3x_1+a^2(d^3-c^3)\cdot x_0x_1^3+c^2(d^3-b^3)\cdot x_0^3x_2\\
                  &+c^2(a^3-d^3)\cdot x_1^3x_2+a^2(b^3-d^3)\cdot x_0x_2^3+b^2(d^3-a^3)\cdot x_1x_2^3\\
                  &+d^2(b^3-c^3)\cdot x_0^3x_3+d^2(c^3-a^3)\cdot x_1^3x_3+d^2(a^3-b^3)\cdot x_2^3x_3\\
                  &+a^2(c^3-b^3)\cdot x_0x_3^3+b^2(a^3-c^3)\cdot x_1x_3^3+c^2(b^3-a^3)\cdot x_2x_3^3. 
   \end{align*}
\endproof
   Let $\eps$ be a primitive root of the unity of degree $3$.
   Then the set $Z$ in Theorem \ref{thm: BMSS} can be written more explicitly as all points of the form
   $$(1:\eps^\alpha:\eps^\beta:\eps^\gamma),\;\mbox{ where }\; \alpha, \beta, \gamma =1,2,3$$
   together with the three coordinate points in $\P^3$. Note that for example the point
   $$(0:0:1:1)$$
   is contained in the $L(\cala_{4}^{4}(3))$, since it is the intersection point of arrangement
   hyperplanes
   $$x_0-x_1=0,\; x_0-\eps x_1=0,\; x_0-\eps^2 x_1=0\;\mbox{ and }\; x_2-x_3=0$$
   but is not an element of $Z$.
\subsection{Unexpected hypersurfaces with multiple general fat points}
   After passing from Definition \ref{def: unexpected curve} to Definition \ref{def: unexpected hypersurface},
   it has been realized that one can allow more general fat points to appear.
\begin{definition}[Unexpected hypersurfaces with multiple points]\label{def: unexpected hypersurface many points}
   We say that a subscheme $Z\subset\P^N$ \emph{admits an unexpected hypersurface} of degree $d$
   with respect to $X=m_1P_1+\ldots +m_sP_s$, where
   $m_1,\ldots,m_s$ are integers and $P_1,\ldots,P_s$ are \textbf{general} points, if
   the fat points scheme $X=m_1P_1+\ldots +m_sP_s$ fails to impose independent conditions on forms
   of degree $d$ vanishing along $Z$, i.e.,
\begin{align*}
   h^0(\P^N,\calo_{\P^N}(d)&\otimes I(Z+\sum{j=1}^{s} m_jP_j))>\\
      &\max\left\{ h^0(\P^N,\calo_{\P^N}(d)\otimes I(Z))-\sum_{j=1}^s\binom{m_j+1}{2},0\right\}.
\end{align*}
\end{definition}
\begin{remark}
   It is relatively easy to construct examples of this kind of behaviour. For instance, let
   $Z$ be an empty set and let $P_1,\ldots,P_7$ be general points in $\P^4$. Then it is well known
   (one of special cases in the Alexander-Hirschowitz classification of special linear systems
   with double base points \cite{AleHir95}) that the scheme $X=2P_1+\ldots+2P_7$
   fails to impose independent conditions on forms of degree $3$. There exists a threefold $T$
   of degree $3$ singular in these $7$ points. However $T$ is singular along the rational quartic curve
   passing through $P_1,\ldots,P_7$. It is much harder to find examples where the points in $X$
   are isolated in the singular locus.
\end{remark}
   Until recently it was not clear that if there exist unexpected hypersurfaces with \emph{isolated} multiple
   general fat points. The first example of this kind has been announced in \cite{Szp19multi}.
   Expectedly, it was constructed with a Fermat-derived configuration of points. More precisely,
   let $Z$ be the union of the six coordinate points in $\P^5$ with the set of all points of the form
   $$(1:\eps^{s_1}:\eps^{s_2}:\eps^{s_3}:\eps^{s_4}:\eps^{s_5})\;\mbox{ with }\; s_1,\ldots,s_5=1,2,3,$$
   where as usually $\eps$ is a primite root of unity of order $3$.
\begin{theorem}\label{thm:P5 with mult 3 and 2}
   Let $R=(a_0:a_1:\ldots:a_5)$ and $S=(b_0:b_1:\ldots:b_5)$ be general points in $\P^5$.
   The set $Z$ as above admits a unique expected quartic
   $4$--fold $Q_{R,S}\subset\P^5$. More precisely $Q_{R,S}$ passes through
   \begin{itemize}
   \item all points in $Z$,
   \item has at $R$ a singularity of order $3$,
   \item has at $S$ a singularity of order $2$.
   \end{itemize}
\end{theorem}

\section{Unexpected curves and Fermat-derived point configurations}
\label{JSsec:5}
   In Section \ref{JSsec: unexpected curves} we have seen how configurations of points dual to lines
   in a Fermat-type arrangement lead to unexpected curves. In the present part we come back
   to Fermat line arrangements and investigate if the configurations of points derived from them
   give rise to unexpected curves. We show that this is indeed the case and somewhat surprisingly
   we discover a new phenomena: No matter how big degree of the Fermat arrangement we consider,
   the multiplicity of an unexpected curve in a general point is always $4$. This is
   in clear opposition to Theorem \ref{thm: Fermat m>=5} where the multiplicity of the unexpected curve
   in a general point grows with $m$. We have no conceptual explanation for this fact at the moment.
\begin{theorem}[Unexpected curves with a point of multiplicity $4$]\label{thm: unexpected curves mult 4}
   Let $Z$ be the configuration of points in $\P^2$ derived from the Fermat arrangement $\cala_3^0(n)$
   for $n\geq 3$.
   Let $P=(a:b:c)$ be a general point in $\P^2$.
   We define the following numbers:
   $$u=\binom{n}{2}-1, \quad v=\binom{n-1}{2}, \quad w=\binom{n+1}{2}.$$
   Then the polynomial
   \begin{align*}
   Q_P(x:y:z)=&-cxy((ub^n+vc^n)(z^n-x^n)+(ua^n+vc^n)(y^n-z^n))\\
        &-bxz((ua^n+vb^n)(y^n-z^n)+(uc^n+vb^n)(x^n-y^n))\\
        &-ayz((ub^n+va^n)(z^n-x^n)+(uc^n+va^n)(x^n-y^n))\\
        &+wa^{n-1}bcx^2(y^n-z^n)+wab^{n-1}cy^2(z^n-x^n)\\
        &+wabc^{n-1}z^2(x^n-y^n) \numberthis\label{eq: mult 4}
   \end{align*}
   \begin{itemize}
   \item vanishes at all points of $Z$,
   \item vanishes to order $4$ at $P$,
   \item defines an unexpected curve of degree $n+2$ for $Z$ with respect to $P$.
   \end{itemize}
\end{theorem}
\proof
   The points in $Z$ are of the form
   $$P_{(\alpha,\beta)}=(1:\eps^\alpha:\eps^\beta )$$
   where $\eps$ is a primitive root of unity of order $n$ and $1\leq\alpha,\beta\leq n$;
   and three coordinate points $P_1=(1:0:0)$, $P_2=(0:1:0)$, $P=3=(0:0:1)$.

   Vanishing of $Q_P$ at the coordinate points is clear, since every summand
   in \eqref{eq: mult 4} vanishes in these points. Similarly, vanishing at points
   $P_{(\alpha,\beta)}$ is guaranteed by vanishing of all summands in \eqref{eq: mult 4}
   in these points.

   The ideal $I(P)$ is generated by
   $$f_1=cx-az,\;\mbox{ and } f_2=cy-bz,$$
   and then the ideal $I(4P)$ is generated by
   $$g_1=f_2^4,\;
     g_2=f_1\cdot f_2^3,\;
     g_3=f_1^2\cdot f_2^2,\;
     g_4=f_1^3\cdot f_2,\;
     g_5=f_1^4.$$
   These generators form in fact a Gr\"obner basis of the ideal $I(4P)$.
   In this basis the polynomial $Q_P$ is presented as follows
   \begin{align*}
   c^4 Q_P &=
         \left((a^4+2ac^3)z-(2a^3c+c^4)\right)\cdot g_1
         +\left(6a^2bc x-(4a^3b+2bc^3)z\right)\cdot g_2+\\
         &+ \left((4ab^3+2ac^3)z-6ab^2cy\right)\cdot g_4
         +\left((2b^3c+c^4)y-(b^4+2bc^3)z\right)\cdot g_5.
   \end{align*}
   Interestingly, the third, "most symmetric", generator of $I(4P)$ is not necessary
   to define the unexpected polynomial $Q_P$.
\endproof

\section{Unexpected surfaces and flats}
   By a flat we mean here a linear subspace of a projective space, i.e.,
   a subscheme defined by linear equations.
   Linear systems with base loci imposed along higher dimensional flats
   have been studied recently by Guardo, Harbourne and Van Tuyl in \cite{GHV13}.
   Their study has been motivated by the containment problem between symbolic
   and ordinary powers of the associated homogeneous ideals. Dumnicki, Harbourne,
   Szemberg and Tutaj-Gasi\'nska in \cite{DHST14} studied linear systems of this
   kind from the positivity point of view. Very recently Migliore, Nagel and Schenck
   in \cite{MigNagSch19}
   initiated a systematic investigation of singular loci of hyperplane arrangements.
   The schemes they consider generalize our concept of Fermat-derived configurations.
   In this spirit we generalize Definition \ref{def: unexpected hypersurface}
   once again.
\begin{definition}[Unexpected hypersurface for flats]
   Let $Z$ be a subscheme of $\P^N$, let $X$ be a general flat in $\P^N$ and let $m$
   be a positive integer.
   We say that $Z$ admits an \emph{unexpected hypersurface} of degree $d$ with respect
   to $mX$ if
   $$h^0(\P^N, \calo_{\P^N}(d)\otimes I(Z+mX))\;>\;
     h^0(\P^N, \calo_{\P^N}(d)\otimes I(Z)) - \HF_X(d),$$
   where $HF_X$ denotes the Hilbert function of the scheme $X$.
\end{definition}
   In other words, there is an unexpected hypersurface if $X$ fails to impose the expected
   number of conditions on linear series.

   Let $c(N,r,m,d)$ be the number of conditions imposed on forms of degree $d$ in $\P^N$
   by vanishing along an $r$--dimensional flat
   to multiplicity at least $m$. These numbers have been computed in \cite[Lemma 2.1]{DHST14}
   and we have the following formula
   \begin{equation*}\label{eq: conditions on fat flat in PN}
      c(N,r,m,d)=\sum_{0\leq i<m}\binom{d-i+r}{r}\binom{N-r-1+i}{i}.
   \end{equation*}
   Specializing to a line in $\P^N$ and rearranging terms we obtain
   \begin{equation*}\label{eq: conditions of fat line in PN}
      c(N,1,m,d)=\frac{m(Nd+2N+m-mN-1)}{N(N-1)}\binom{N+m-2}{m}.
   \end{equation*}
   And finally specializing to a line in $\P^3$ we get
   \begin{equation*}\label{eq: conditions of fat line in P3}
      c(m,d)=c(3,1,m,d)=\binom{m+1}{2}(t+1)-2\binom{m+1}{3}.
   \end{equation*}
   It is well known that a single line with arbitrary multiplicity imposes
   independent conditions on forms of any degree. This is in parallel with
   conditions imposed by a single point. A ground breaking discovery of
   Cook II, Harbourne, Migliore and Nagel in \cite{CHMN} was that
   it may happen that a fat point imposes \emph{dependent}
   conditions in a \emph{noncomplete} linear system (of hypersurfaces
   vanishing along $Z$). It is clear that a general point must be taken with
   multiplicity at least $2$ in order to exhibit this kind of behavior.
   In fact Akesseh noticed in \cite{Ake17} that a general point of multiplicity $2$
   can impose dependent conditions on a noncomplete linear series only
   if the characteristic of the ground field is $2$. Farnik, Galuppi, Sodomaco
   and Trok extended this picture showing that the only unexpected curve
   with a point of multiplicity $3$ is the $B_3$ quartic given in \eqref{eq: B3 quartic}.

   We show now that, surprisingly, it may happen that a single general line fails
   to impose independent conditions on forms with a base locus consisting of a Fermat-derived
   configuration of lines in $\P^3$. Our observation is experimental based
   on Singular computations.

   Let $Z\subset \P^3$ be the union of those lines derived from Fermat-type arrangement $\cala_4^0(3)$,
   where at least three of arrangement planes intersect. It has been computed in \cite[Section 3.1]{MalSzp18}
   that there are $42$ mutually distinct (but not disjoint) lines in $Z$. By \cite[Lemma 4.1]{MalSzp18} the ideal $I(Z)\subset R=\C[x,y,z,w]$
   is generated by
   $$(x^3-y^3)(z^3-w^3)xy,\; (x^3-y^3)(z^3-w^3)zw,$$
   $$(x^3-z^3)(y^3-w^3)xz,\; (x^3-z^3)(y^3-w^3)yw,$$
   $$(x^3-w^3)(y^3-z^3)xw,\; (x^3-w^3)(y^3-z^3)yz.$$
   It is easy to check, that the minimal free resolution of $I(Z)$ has the form
   $$I(Z)\leftarrow R(-8)^6\leftarrow R(-9)^4\oplus R(-12)\leftarrow 0.$$
   In particular the dimension of the space of forms of degree $9$ in $I(Z)$ is exactly $20$.
   Since a reduced line imposes $c(1,d)=d+1$ conditions on forms of degree $10$,
   we expect that for a general line $L\subset\P^3$ the space of
   forms of degree $9$ in $I(Z+L)$ will be $10$. However, according to Singular, this space
   has dimension $12$, which is unexpected. To be more precise we run computations with
   a \emph{random} line rather than a general line but nevertheless the indications that
   $12$ is indeed the dimension are rather strong. We were not able to verify the computations
   for a general line because it never stopped. It would be, of course, desirable
   to have a theoretical explanation for this phenomena and a more general phenomena, which
   we state now.
\begin{observation}
   Let $Z$ be the set of lines as above. Let $m$ be a positive integer and let $L$ be a general (random)
   line in $\P^3$. Then $Z$ admits unexpected surfaces of degree $m+8$ with respect to the scheme $mL$
   for all $m\geq 1$.
\end{observation}
\begin{remark}
   Unfortunately, in opposition to results presented so far, we were not able to write down explicitly
   equations of the unexpected surfaces. In fact, Singular computations break, if the line is supposed
   to be defined by two general linear equations of the form
   $$ax+by+cz+dw=0.$$
   We hope to come back to this problem in the very near future.
\end{remark}

\section{Concluding remarks}
   In this note I made an attempt to introduce in a systematic way Fermat-type arrangements of hyperplanes
   in projective spaces and configurations of flats (including points) derived as intersections
   among those hyperplanes. I showed that these arrangements are always free since they are
   reflection arrangements. I also explained how they constitute an ample group of examples in the theory
   of unexpected hypersurfaces. Since they play also an important role in the containment problem
   for powers of ideals in commutative algebra, it can be expected that there might more areas
   of algebra, geometry and combinatorics where Fermat-type arrangements appear as interesting source
   of examples. This work can be regarded as an invitation to study these interesting objects
   in greater detail and to consider them as a rich testing field for various statements and ideas.

\paragraph*{Acknowledgement.}
   This research was partially supported by National Science Centre, Poland, grant
   2018/02/X/ST1/00519. I would like to thank Tomasz Szemberg for suggesting to me to write up
   this survey. I thank also Marcin Dumnicki and Piotr Pokora for helpful remarks
   on the first draft of the note.



\end{document}